# 2-EDGE DISTANCE-BALANCED GRAPHS

**Zohreh Aliannejadi[1], Mehdi Alaeiyan[2], Alireza Gilani[1] & Jafar Asadpour[1]**

[1]Department of Mathematics, South Tehran Branch
Islamic Azad University, Tehran, Iran
[2]Department of Mathematics, Iran University of Science and Technology
Narmak, Tehran 16844, Iran
[*]Corresponding author: email address alaeiyan@iust.ac.ir

**Abstract**
In a graph $A$, for each two arbitrary vertices $g, h$ with $d(g,v) = 2$, $|M^A_{g_2h}| = m^A_{g_2h}$ is introduced the number of edges of $A$ that are closer to $g$ than to $h$. We say $A$ is a 2-edge distance-balanced graph if we have $m^A_{g_2h} = m^A_{h_2g}$. In this article, we verify the concept of these graphs and present a method to recognize k-edge distance-balanced graphs for $k = 2,3$ using existence of either even or odd cycles. Moreover, we investigate situations under which the cartesian and lexicographic products lead to 2-edge distance-balanced graphs. In some subdivision-related graphs 2-edge distance-balanced property is verified.
Mathematics Subject Classification: 05C12, 05C25

**Keywords:** Cartesian product; Lexicographic product; Subdivision graphs; Total distance; 2-edge distance-balanced graphs;

## 1 INTRODUCTION

The notion of graph is a pivotal tool to make use of the modeling of the phenomena and it is taken into consideration in many studies in a recent decades. One of the optimal uses of graphs theory is to classify graphs based on discriminating quality. This phenomenon can be best observed in distance-balanced graphs has been determined by [10]. Also, it is investigated in some papers, we refer the reader to ([1],[2],[3],[5],[8],[11]-[15]) and references therein.

Let $A$ be a connected, finite and undirected graph throughout of this paper, in which $V(A)$ is its vertex set and $E(A)$ is its edge set. In a graph $A$, the distance between each pair of vertices $g, h \in V(A)$ is introduced the number of edges in the least distance joining them and it is indicated by $d_A(g,h)$ (see [3, 15]. For every two desired edges $f = gh, \acute{f} = \acute{g}\acute{h}$, the distance between $f$ and $\acute{f}$ is introduced via:
$d_A(f,\acute{f}) = \min\{d_A(g,\acute{f}), d_A(h,\acute{f})\}$
$= \min\{d_A(g,\acute{g}), d_A(g,\acute{h}), d_A(h,\acute{g}), d_A(h,\acute{h})\}$.
Set $m_g(f) = |M_g(f)| = |\{\acute{f} \in E(A) | d_A(g,\acute{f}) < d_A(h,\acute{f})\}|$

$m_v(f) = |M_h(f)| = |\{\acute{f} \in E(A) | d_A(h,\acute{f}) < d_A(g,\acute{f})\}|$
and $m_0(f) = |M_0(f)| = |\{\acute{f} \in E(G) | d_A(g,\acute{f}) = d_A(h,\acute{f})\}|$.

Presume that $f = gh \in E(A)$. For every two integers $i, j$, we consider:
$\acute{D}^i_j(e) = \{\acute{f} \in E(A) | d_A(\acute{f},g) = i. d_A(\acute{f},h) = j\}$.
A "distance partition" of $E(A)$ is concluded by sets $\acute{D}^i_i(f)$ due to the edge $f = gh$. Merely the sets $\acute{D}^{i-1}_i(f), \acute{D}^i_i(f)$ and, $\acute{D}^i_{i-1}(f)$, for every $(1 \leq i \leq d)$ might be nonempty according to the triangle inequality ($d$ is the diameter of the graph $A$). As well as $\acute{D}^0_0(e) = \phi$.

For an edge $f = gh$ of $A$ we denote $n^A_g(f) = |W^A_{g,h}| = |\{a \in V(A) | d_A(a,g) < d_A(a,h)\}|$.
Analogously, we would define $n^A_h(f) = |W^A_{h,g}|$.
Graph $A$ is named distance-balanced (DB) whenever for an edge $f = gh$ of $A$ we have $n^A_g(f) = n^A_h(f)$ [10].

Graph $A$ is called edge-distance-balanced (EDB) whenever for an edge $f = gh$ of $A$ there is a positive integer $\gamma$, so that $m^A_g(f) = m^A_h(f) = \gamma$ [16].

In Section 2, it is introduced an extended version of the concept of edge distance-balanced, that is called 2-edge distanced-balanced. We study 2-edge distance-balanced graphs in the framework of cartesian and lexicographic products of two connected graphs in Section 3. Finally, in Section 4, it is investigated 2-edge distance-balanced property in some subdivision-related graphs.

**2 Specification of k-edge distance-balanced graphs ($k = 2, 3$)**

In this segment, we would complete results of n-distance-balanced graphs that it has already expressed by [6] for k-edge distance-balanced graphs for $k = 2,3$ and present a method for classification such graphs.

**Definition 21** Let $A$ be a graph. It is denominated 2-edge distance-balanced (briefly 2 – EDB) if and only if for each two of vertices $g, h \in V(A)$ with $d(g, h) = 2$, it holds that $|M^A_{g\underline{2}h}| = |M^A_{h\underline{2}g}|$, in which
$$m^A_{g\underline{2}h} = |M^A_{g\underline{2}h}| = |\{f \in E(A) | d(f, g) < d(f, h)\}|.$$
similarly, $m^A_{h\underline{2}g} = |M^A_{h\underline{2}g}| = |\{f \in E(A) | d(f, h) < d(f, g)\}|$.
Also, consider the concept $m^A_{0g\underline{2}h} = |M^A_{0g\underline{2}h}| = \{f \in E(A) | d(f, h) = d(f, g)\}$.

The notion $g\underline{2}h$ is equivalent a path with length 2. Obviously, a distance partition for $E(A)$ is formed by $M^A_{g\underline{2}h}$, $M^A_{h\underline{2}g}$ and $M^A_{0\,g\underline{2}h}$.

**Example 22** Star graphs $S_k$ (with $k > 1$), Wheel graphs $W_n$ (with $n > 4$), Frienship graphs $F_n$ (with $n > 1$) are 2-EDB graphs. Complete bipartite graph $K_{m.n}$ is 2-EDB but without edge distance-balanced property.

**Definition 23** The total distance $D^A(g)$ of $g$, for a vertex $g$ of $A$ is denoted
$D^A(g) = \sum_{f \in E(A)} d^A(g, f) = d(g, E(A))$.
Whenever $A$ is explicit from the text we will write $d(g, f)$ and $D(g)$ in lieu of $d^A(g, f)$ and $D^A(g)$, respectively.

**Definition 24** The property $(\Delta_n)$ in a connected graph $A$, for vertices $g, h \in V(A)$ with $d(g, h) = n$ and $f \in M^A_{g\underline{n}h}$ is introduced as the shortest way between $f$ and $h$ that is, $(W_{fh})$ so that it does not include the shortest way between $f$ and $h$ that is, $(W_{fh})$.

**Theorem 1** Presume that $A$ is a connected graph. Then under one of the below conditions, $A$ is 2-EDB if and only if for each two of vertices $g, h \in V(A)$ with $d(g, h) = 2$, $D^A(g) = D^A(h)$ holds.
(i) $A$ does not have any odd cycle, that is, $A$ is bipartite.
(ii) $A$ does not have any even cycle, but has the property $(\Delta_2)$.

**Proof** Since condition (i) holds, we assume that $g, h$ are two vertices in $A$, in which $d(g, h) = 2$. Then $D^A(g) = D^A(h)$ can be shown as follows
$\sum_{f \in M^A_{g\underline{2}h}} d(g, f) + \sum_{e \in M^A_{h\underline{2}g}} d(g, f) +$
$\sum_{f \in M^G_{0\,g\underline{2}h}} d(g, f) =$
$\sum_{f \in M^A_{g\underline{2}h}} d(h, f) + \sum_{f \in M^A_{h\underline{2}g}} d(h, f) +$
$\sum_{f \in M^A_{0\,g\underline{2}h}} d(h, f).$
which concludes that
$$\sum_{f \in M^A_{g\underline{2}h}} (d(g, f) - d(h, f)) = \sum_{f \in M^A_{h\underline{2}g}} (d(h, f) - d(g, f)). \qquad (1)$$
For any $f \in M^A_{g\underline{2}h}$ we have $d(g, f) - d(h, f) = -2$. To prove it, we consider that $f \in M^A_{g\underline{2}h}$ and $d(f, g) = l(P)$ so that $P$ is the shortest path between $f$ and $g$.
Hence,
$l(P) = d(f, g) < d(f, h) \leq l(P) + 2$.
It implies that $d(f, h) = l(P) + 1$ or $l(p) + 2$. So long as $d(f, h) = l(P) + 1$, it can be inferred that there exists an odd cycle having length $2l(P) + 1$ or $2l(P) + 3$. Hence, it contradicts and the claim follows. In the same way, it attains $d(h, f) - d(g, f) = -2$ for each $f \in M^A_{h\underline{2}g}$. Now by (1) we obtain
$\sum_{f \in M^A_{g\underline{2}h}} (-2) = \sum_{f \in M^A_{h\underline{2}g}} (-2).$
Thus, it is right if and only if $A$ is 2-EDB. Hence, under condition (i) the proof is completed. For next condition, we consider that $A$ does not have any even cycle and has the property $(\Delta_2)$. Since $d(f, g) = l(P)$, for each edge $f \in M^A_{g\underline{2}h}$ then we obtain $d(f, h) \in \{l(P) + 1. l(P) + 2\}$. Since $d(f, h) = l(P) + 2$, then a path $R$ between f and $h$ with length $l(P) + 2$ exists. It is seen that there is an even cycle with length $2l(P) + 4$ or $2l(P) + 2$ in $A$ that is a contradiction. Thus,
$\sum_{f \in M^A_{g\underline{2}h}} (d(g, f) - d(h, f)) =$
$\sum_{f \in M^A_{h\underline{2}g}} (d(h, f) - d(g, f)).$
then we attain $\sum_{f \in M^A_{g\underline{2}h}} (-1) = \sum_{f \in M^A_{h\underline{2}g}} (-1)$.
The proof ends. □

Theorem 1 matters to classify a group of the edge distance-balanced graphs with specified properties.

**Corollary 25** Every bipartite and 2-EDB graph is edge distance-balanced.

In the next theorem, using the notion of total distance, we present an analogous condition and a procedure to recognize edge distance-balanced graphs with $d(g,h) = 3$ which are called 3-EDB graphs.

**Theorem 2** If $A$ be a connected graph having the property $(\Delta_3)$ without any even cycle, then $A$ is 3-EDB if and only if for each two of vertices $g, h \in V(A)$ with $d(g,h) = 3$, it holds $D^A(g) = D^A(h)$.

**Proof** Based on the proof of Theorem 1, we assert that
$$\sum_{f \in M^A_{g_3h}} (d(g,f) - d(h,f)) = \sum_{f \in M^A_{h_3g}} (d(h,f) - d(g,f))$$
and it follows $\sum_{f \in M^A_{g_3h}} (-2) = \sum_{f \in M^A_{h_3g}} (-2)$.

Let $f \in M^A_{g_3h}$ and $d(f,g) = l(R)$ such that $R$ is the shortest path between $f$ and g. If $l(R) = d(g,f) < d(h,f) \le l(R) + 3$, then $d(h,f)$ is $l(R)+1$ or $l(R)+2$ and or $l(R)+3$. Suppose that $m$ and $n$ are in the path between $g$ and $h$. We call $l(Q) = d(f,m)$, $l(S) = d(f,n)$ and $l(P) = d(f,h)$. Since $d(h,f) = l(R) + 1$, then concerning the property $(\Delta_3)$ and with paths $Q, P$ and $S$ we obtain an even cycle with $2l(R)$, $2l(R) + 2$ or $2l(R) + 4$, respectively. Hence it contradicts. Analogously, since $d(f,h) = l(R) + 3$, then applying the paths $Q, S$ or $P$ an even cycle is observed with length $2l(R) + 2, 2l(R) + 4$ or $2l(R) + 6$, respectively which makes a contradiction, as well as. So the assertion is proved. Therefore, the result follows. □

## 3 2-edge distance-balanced graphs and product graphs

We would now investigate situations in which the *Cartesian product* leads to 2-EDB graphs. We mention that such product graphs, formed by graphs $A$ and $B$, its vertex set is $(A \square B) = V(A) \times V(B)$. Consider that $(a_1, b_1)$ and $(a_2, b_2)$ are distinct vertices in $V(A \square B)$. In the Cartesian product $A \square B$, if two vertices $(a_1, b_1)$ and $(a_2, b_2)$ are coincident in one coordinate and adjacent in the other coordinate, then they are adjacent, that is, $a_1 = a_2$ and $b_1 b_2 \in E(B)$, or $b_1 = b_2$ and $a_1 a_2 \in E(A)$. Obviously, for vertices we have:
$$d_{A \square B}((a_1, b_1), (a_2, b_2)) = d_A(a_1, a_2) + d_B(b_1, b_2).$$

For edges we attain:
$$d_{A \square B}((a,b)(a_1,b_1),(á,b́)(á_1,b́_1)) =$$
$$\min\{d_{A \square B}((a,b),(á,b́)), d_{A \square B}((a,b),(á_1,b́_1)), d_{A \square B}((a_1,b_1),(á,b́)), d_{A \square B}((a_1,b_1),(á_1,b́_1))\} =$$
$$\min\{d_A(a,á) + d_B(b,b́), d_A(a,á_1) + d_B(b,b́_1), d_A(a_1,á) + d_B(b_1,b́), d_A(a_1,á_1) + d_B(b_1,b́_1)\}.$$

**Theorem 3** If $A$ and $B$ are graphs, then $A \square B$ is 2-EDB if and only if both $A$ and $B$ are 2-EDB and 2-DB.

**Proof** Consider that $a_1, a_2$ are two adjacent vertices in $A$ with $d(a_1, a_2) = 2$ and also $b_1, b_2$ are two vertices in $B$ with $d(b_1, b_2) = 2$. Let $(a_1, b_1)$, $(a_2, b_1)$ and $(a_1, b_2) \in V(A \square B)$. Then it is clearly seen that
$$d((a_1,b_1),(a_2,b_1)) = 2, d((a_1,b_1),(a_1,b_2)) = 2 \text{ and } d((a_2,b_1),(a_1,b_2)) = 2.$$
We see that

$$M_{(a_1,b_1)_2(a_2,b_2)} = \{(a,b)(á,b́) \in E(A \square B) | aá \in E(A), b = b́ \text{ or } bb́ \in E(B), a = á.$$
$$\min\{d_A(a,a_1), d_A(á,a_1)\} < \min\{d_A(a,a_2), d_A(á,a_2)\}\}. \quad (2)$$
and hence applying (2) and by [16. Theorem 2.1) we conclude that:
$$m^{A \square B}_{(a_1,b_1)_2(a_2,b_2)} = m^A_{a_1 2 a_2} \cdot |V(B)| + n^A_{a_1 2 a_2} \cdot |E(B)|. \quad (3)$$
Similar to this process we obtain
$$m^{A \square B}_{(a_2,a_2)_2(a_1,b_1)} = m^A_{a_1 2 a_2} \cdot |V(B)| + n^A_{a_1 2 a_2} \cdot |E(B)|, \quad (4)$$

$$m^{A\square B}_{(a_1,b_1)\underline{2}(a_1,b_2)} = m^B_{b_1\underline{2}b_2} \cdot |V(A)| + n^B_{b_1\underline{2}b_2} \cdot |E(A)|, \quad (5)$$

$$m^{A\square B}_{(a_1,b_2)\underline{2}(a_1,b_1)} = m^B_{b_2\underline{2}b_1} \cdot |V(A)| + n^B_{b_2\underline{2}b_1} \cdot |E(A)|. \quad (6)$$

Now, suppose that $A$ and $B$ are 2-EDB and 2-DB graphs. By (3) and (4) we have:
$$m^A_{a_1\underline{2}a_2} \cdot |V(B)| + n^A_{a_1\underline{2}a_2} \cdot |E(B)|$$
$$= m^A_{a_2\underline{2}a_1} \cdot |V(B)| + n^A_{a_2\underline{2}a_1} \cdot |E(B)|.$$

therefore
$$m^{A\square B}_{(a_1,b_1)\underline{2}(a_2,b_2)} = m^{A\square B}_{(a_2,b_2)\underline{2}(a_1,b_1)}.$$

By analogy, using (5) and (6), and we know that B is 2-EDB and 2-DB, we conclude that
$$m^{A\square B}_{(a_1,b_1)\underline{2}(a_1,b_2)} = m^{A\square B}_{(a_1,b_2)\underline{2}(a_1,b_1)}.$$
and therefore $A\square B$ is 2-EDB.

For converse, consider that $A\square B$ is 2-EDB, then using (3) and (4) we observe that

$$m^{A\square B}_{(a_1,b_1)\underline{2}(a_2,b_2)} = m^{A\square B}_{(a_2,b_2)\underline{2}(a_1,b_1)} \implies$$
$$m^A_{a_1\underline{2}a_2} \cdot |V(B)| + n^A_{a_1\underline{2}a_2} \cdot |E(B)|$$
$$= m^A_{a_2\underline{2}a_1} \cdot |V(B)| + n^A_{a_2\underline{2}a_1} \cdot |E(B)|.$$

therefore $A$ is 2-EDB and 2-DB. Similarly (5) and (6) yield that $B$ is 2-EDB and 2-DB. The proof is completed. □

Here, we would define the lexicographic product graphs. The *lexicographic product* $A[B]$ of two graphs $A$ and $B$ is the graph that $V(A[B]) = V(A) \times V(B)$ is its vertex set and two vertices $(a_1, b_1)$, $(a_2, b_2)$ are adjacent if $a_1 a_2 \in E(A)$ or if $a_1 = a_2$ and $b_1 b_2 \in E(B)$ [9, p.22]. Since $A$ is a connected non-trivial graph, then it is easily seen for vertices that

$$d_{A[B]}((a_1,b_1),(a_2,b_2)) = \begin{cases} d_A(a_1,a_2) & if\, a_1 \neq a_2 \\ \min\{2, d_B(b_1,b_2)\} & if\, a_1 = a_2. \end{cases}$$

And for edges we have:
$$d_{A[B]}((a,b)(a_1,b_1),(á,b́)(á_1,b́_1)) = \min \begin{cases} d_A(a,á) & if\, a \neq á, & \min\{2, d_B(b,b́)\} & if\, a = á \\ d_A(a,á_1) & if\, a \neq á_1, & \min\{2, d_B(b,b́_1)\} & if\, a = á_1 \\ d_A(a_1,á) & if\, a_1 \neq á, & \min\{2, d_B(b_1,b́)\} & if\, a_1 = á \\ d_A(a_1,á_1) & if\, a_1 \neq á_1, & \min\{2, d_B(b_1,b́_1)\} & if\, a_1 = á_1 \end{cases}.$$

In the following theorem, it is proved that the lexicographic product $A[B]$ of graphs $A$ and $B$ is 2-EDB if and only if $A$ is 2-EDB and $B$ is locally regular. At first, we state the definition of locally regular graphs.

**Definition 31** The graph $A$ is locally regular according to $n$ (briefly n-locally regular) if it holds
$\forall g, h \in V(A), \quad d(g,h) = n \Rightarrow deg(g) = deg(h)$ [6].

**Theorem 4** Presume that $A$ and $B$ are graphs. Then the graph $A[B]$ is 2-EDB if and only if $A$ is 2-EDB and $B$ is locally regular.

*Proof* For the beginning, let the graph $A[B]$ be 2-EDB and $(a_1, b_1), (a_1, b_2) \in V(A[B])$, $b_1 \neq b_2$ with $d((a_1, b_1), (a_1, b_2)) = 2$. Then $b_1 b_2 \notin E(A)$. For edge $(á_1, ḣ_1)(á_2, ḃ_2) \in E(A[B])$ according to the definition of distance between an edge and a vertex, we obtain the following consequence

$$d_{A[B]}((á_1, ḃ_1)(á_2, ḃ_2), (a_1, b_1)) < d_{A[B]}((á_1, ḃ_1)(á_2, ḃ_2), (a_1, b_2)) \Rightarrow$$
$$\min\{d_{A[B]}((á_1, ḃ_1), (a_1, b_1)), d_{A[B]}((á_2, ḃ_2), (a_1, b_1))\} <$$
$$\min\{d_{A[B]}((á_1, ḃ_1), (a_1, b_2)), d_{A[B]}((á_2, ḃ_2), (a_1, b_2))\} \Rightarrow$$
$$\min\{d_{A[B]}((á_1, ḃ_1), (a_1, b_2)), d_{A[B]}((á_2, ḃ_2), (a_1, b_2))\} \Rightarrow$$

$$\begin{cases} á_1 = a_1, & ḃ_1 b_1 \in E(B), ḃ_1 b_2 \notin E(B) \\ \quad \text{or} \\ á_2 = a_1, & ḃ_2 b_1 \in E(B), ḃ_2 b_2 \notin E(B). \end{cases} \quad (7)$$

Similarly,
$$d_{A[B]}((á_1, ḃ_1)(á_2, ḃ_2), (a_1, b_2)) < d_{A[B]}((á_1, ḃ_1)(á_2, ḃ_2), (a_1, b_1)) \Rightarrow$$
$$\min\{d_{A[B]}((á_1, ḃ_1), (a_1, b_2)), d_{A[B]}((á_2, ḃ_2), (a_1, b_2))\} <$$
$$\min\{d_{A[B]}((á_1, ḃ_1), (a_1, b_1)), d_{A[B]}((á_2, ḃ_2), (a_1, b_1))\} \Rightarrow$$

$$\begin{cases} á_1 = a_1, & ḃ_1 b_2 \in E(B), ḃ_1 b_1 \notin E(B) \\ \quad \text{or} \\ á_2 = a_1, & ḃ_2 b_2 \in E(B), ḃ_2 b_1 \notin E(B). \end{cases} \quad (8)$$

Since $A[B]$ is 2-EDB by equalities (7) and (8) it follows that,
$$|\{ḃ_1 | ḃ_1 b_1 \in E(A)\}| = |\{ḃ_1 | ḃ_1 b_2 \in E(A)\}|,$$
and
$$|\{ḃ_2 | ḃ_2 b_1 \in E(A)\}| = |\{ḃ_2 | ḃ_2 b_2 \in E(A)\}|.$$

If $b_1 b_2 \notin E(A)$, then the above implication follows that any non-adjacent vertices of $B$ have the same degree. Hence $B$ is a locally regular graph.

Now suppose that $(a_1, b_1), (a_2, b_2) \in V(A[B])$ and $(á_1, ḃ_1)(á_2, ḃ_2) \in E(A[B])$, in which $a_1 \neq a_2$ and $d((a_1, b_1), (a_2, b_2)) = 2$. Then $d(a_1, a_2) = 2$ and it implies
$(á_1, ḃ_1)(á_2, ḃ_2) \in M^{A[B]}_{(a_1, b_1)\underline{2}(a_2, b_2)} \Leftrightarrow$ \hfill (9)
$$d_{A[B]}((á_1, ḃ_1)(á_2, ḃ_2), (a_1, b_1)) < d_{A[B]}((á_1, ḃ_1)(á_2, ḃ_2), (a_2, b_2)) \Leftrightarrow$$
$$\min\{d_{A[B]}((á_1, ḃ_1), (a_1, b_1)), d_{A[B]}((á_2, ḃ_2), (a_1, b_1))\} <$$
$$\min\{d_{A[B]}((á_1, ḃ_1), (a_2, b_2)), d_{A[B]}((á_2, ḃ_2), (a_2, b_2))\} \Leftrightarrow$$
$$\begin{cases} a_1 \neq á_1 \Leftrightarrow á_2 \notin \{a_1, a_2\} \\ a_1 \neq á_2 \Leftrightarrow á_1 \notin \{a_1, a_2\} \end{cases} \Leftrightarrow á_1 á_2 \in M^A_{a_1 \underline{2} a_2}.$$

Similarly,
$(á_1, ḃ_1)(á_2, ḃ_2) \in M^{A[B]}_{(a_2, b_2)\underline{2}(a_1, b_1)} \Leftrightarrow$ \hfill (10)
$$\begin{cases} a_2 \neq á_1 \Leftrightarrow á_1 \notin \{a_1, a_2\} \\ a_2 \neq á_2 \Leftrightarrow á_2 \notin \{a_1, a_2\} \end{cases} \Leftrightarrow á_1 á_2 \in M^A_{a_2 \underline{2} a_1}.$$

By (9) and (10) and we know that $A[B]$ is 2-EDB, it implies that $A$ is a 2-EDB graph. For converse, consider that $A$ is 2-EDB and $B$ is locally regular. If $(a_1, b_1), (a_2, b_2) \in V(A[B])$, $a_1 \neq a_2$ and $((a_1, b_1), (a_2, b_2)) = 2$. Then $d(a_1, a_2) = 2$. Therefore, we have both relations (9) and (10) which are inferred that $A[B]$ is. 2-EDB. Thus, $A[B]$ is 2-EDB if and only if $A$ is 2-EDB. Consider that $(a_1,b_1),(a_1,b_2) \in V(A[B])$, $b_1 \neq b_2$ and $((a_1,b_1),(a_1,b_2)) = 2$. Then $b_1 b_2 \notin E(A)$ and

$$(\acute{a}_1,\acute{b}_1)(\acute{a}_2,\acute{b}_2) \in M^{A[B]}_{(a_1,b_1)\underline{2}(a_1,b_2)} \Leftrightarrow$$

$$\begin{cases} \acute{a}_1 = a_1, \ b_1\acute{b}_1 \in E(B) \ \text{and} \ b_2\acute{b}_1 \notin E(B) \\ \acute{a}_2 = a_1, \ b_1\acute{b}_2 \in E(B) \ \text{and} \ b_2\acute{b}_2 \notin E(B). \end{cases}$$

also

$$(\acute{a}_1,\acute{b}_1)(\acute{a}_2,\acute{b}_2) \in M^{A[B]}_{(a_1,b_2)\underline{2}(a_1,b_1)} \Leftrightarrow$$

$$\begin{cases} \acute{a}_1 = a_1, \ b_2\acute{b}_1 \in E(B) \ \text{and} \ b_1\acute{b}_1 \notin E(B) \\ \acute{a}_2 = a_1, \ b_2\acute{b}_2 \in E(B) \ \text{and} \ b_1\acute{b}_2 \notin E(B). \end{cases}$$

Therefore,

$$M^{A[B]}_{(a_1,b_1)\underline{2}(a_1,b_2)} = \{(a_1,\acute{b}_1)(a_1,\acute{b}_2) | b_1\acute{b}_1, b_1\acute{b}_2 \in E(B), \ b_2\acute{b}_1, b_2\acute{b}_2 \notin E(B)\}.$$

$$M^{A[B]}_{(a_1,b_2)\underline{2}(a_1,b_1)} = \{(a_1,\acute{b}_1)(a_1,\acute{b}_2) | b_2\acute{b}_1, b_2\acute{b}_2 \in E(B), \ b_1\acute{b}_1, b_1\acute{b}_2 \notin E(B)\}.$$

and the fact that $B$ is locally regular follows that $A[B]$ is 2-EDB and the proof is completed. □

## 42-EDB property in some subdivision-related graphs

In the following, we are going to investigate 2-EDB property in some subdivision-related graphs. Now we need to introduce two subdivision-related graphs which are called $R(A)$ and $S(A)$ and have been introduced by [4, 17].

$S(A)$ is constructed from $A$, such that an extra vertex inserts in each edge of $A$. In other words, a path of length 2 replaces any edge of $A$. (See Fig. 1, 2 for an example).

$R(A)$ is constructed from $A$ by adding a new vertex corresponding to any edge of $A$, then connecting any new vertex to the end vertices of the corresponding edge. Replacing any edge of $A$ by a triangle is another way to express $R(A)$. (See Fig. 1,3 for an example).

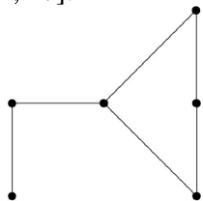
Figure 1: Graph $A$

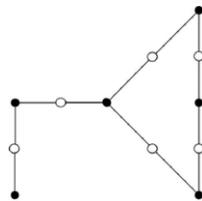
Figure 2: Graph $S(A)$

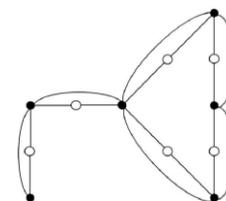
Figure 3: Graph $R(A)$

**Corollary 41** A graph $A$ is EDB if and only if $S(A)$ is 2-EDB.

**Theorem 5** Consider that $A$ is a connected graph. Then $R(A)$ is 2-EDB if and only if $A$ is a path with $|E(A)| = 2$.

*Proof* Suppose that $A$ is a graph, in which $|E(A)| = 2$. So it is easily seen $R(A)$ is a Friendship graph and consequently 2-EDB. For converse, by assumption contrary we consider that $R(A)$ is 2-EDB, such that $A$ be a graph with $|E(A)| > 2$. Thus, there is at least a pair of vertices $g$ and $h$ in $A$ with $d(g,h) = 2$. Suppose that there is a vertex $m$ between $g$ and $h$. Based on the definition of $(A)$, there is a new vertex $x$ corresponding edge $gm$ in $R(A)$. Then, edge $gm$ is one of the edges of constructed triangle in $(A)$. It is seen that $deg(g) = 2$. Similarly, there is a new vertex $y$ corresponding edge $mh$. Then, $mh$ is one of the edges of another triangle in $(A)$. Assume that $m$ and $n$ are another two vertices in $A$ with $d(m,n) = 2$, such that $h$ is between $m$ and $n$. In $(A)$, there is a new vertex $z$ corresponding edge $hn$. Hence, $deg(h) \geq 3$ and $m^{R(A)}_{g\underline{2}h} \neq m^{R(A)}_{h\underline{2}g}$ and it is contradiction. Thus, it is completed. □

**Theorem 6** If $A$ be a nontrivial connected graph with a pendant, then $S(A)$ is not 2-EDB.

*Proof* Let $g, h \in V(A)$ with $d(g,h) = 2$ and $u$ be a pendant vertex. Let $m$ be the new vertex, where $g$ and $m$ and also two vertices $h$ and $m$ are adjacent in $(A)$. Then $m^{S(A)}_{g\underline{2}h} = 1$ and $m^{S(A)}_{h\underline{2}g} \geq 2$. Therefore, $m^{S(A)}_{g\underline{2}h} \neq m^{S(A)}_{h\underline{2}g}$, proving the result.